\documentclass[my_ref]{my2_birkmult}

\newtheorem{thm}{Theorem}[section]

\theoremstyle{definition}
\newtheorem{defn}[thm]{Definition}

\theoremstyle{remark}

\newtheorem{ex}{Example}

\numberwithin{equation}{section}

\usepackage{cite}

\renewcommand{\subjclassname}{Mathematics Subject Classification (2010)}

\begin{document}

\title{Fractional variational calculus of variable order}

\author{T. Odzijewicz}

\address{Center for Research and Development in Mathematics and Applications\\
Department of Mathematics, University of Aveiro, 3810-193 Aveiro, Portugal}

\email{tatianao@ua.pt}

\thanks{Part of the first author's Ph.D., which is carried out
at the University of Aveiro under the
\emph{Doctoral Programme in Mathematics and Applications}
of Universities of Aveiro and Minho.}


\author{A. B. Malinowska}

\address{Faculty of Computer Science, Bia{\l}ystok University of Technology\\
15-351 Bia\l ystok, Poland}

\email{a.malinowska@pb.edu.pl}


\author{D. F. M. Torres}

\address{Center for Research and Development in Mathematics and Applications\\
Department of Mathematics, University of Aveiro, 3810-193 Aveiro, Portugal}

\email{delfim@ua.pt}




\date{}

\dedicatory{To Professor Stefan Samko on the occasion of his 70th birthday}


\begin{abstract}
We study the fundamental problem of the calculus of variations
with variable order fractional operators.
Fractional integrals are considered
in the sense of Riemann--Liouville
while derivatives are of Caputo type.
\end{abstract}

\maketitle


\section{Introduction}

Fractional calculus is a discipline that studies integrals
and derivatives of non-integer (real or complex) order
\cite{book:Kilbas,book:Klimek,book:Podlubny,Samko:book}.
The subject is nowadays very active due to its many applications
in mechanics, chemistry, biology, economics, and control theory \cite{TM:rec}.
In 1993, Samko and Ross proposed an interesting generalization
to fractional calculus \cite{Samko:Ross2} (see also Samko's paper
\cite{Samko} of 1995). They introduced the study
of fractional integration and differentiation
when the order is not a constant but a function.
More precisely, they considered an extension
of Riemann--Liouville and Fourier definitions
\cite{Samko:Ross1,Samko:Ross2,Samko}. Afterwards, several works
were dedicated to variable order fractional operators, their
applications and interpretations (see, \textrm{e.g.}, \cite{alexandre,Coimbra,Lorentzo}).
In particular, Samko's variable order fractional calculus turns out to be
very useful in mechanics and in the theory of viscous flows
\cite{Coimbra,Diaz:Coimbra,Lorentzo,Ramirez:Coimbra1,Ramirez:Coimbra2,Pedro}.
Indeed, many physical processes exhibit fractional-order behavior that may vary with time or space
\cite{Lorentzo}. The paper \cite{Coimbra} is devoted to the study of a variable-order fractional
differential equation that characterizes some problems in the theory of viscoelasticity.
In \cite{Diaz:Coimbra} the authors analyze the dynamics and control of a nonlinear
variable viscoelasticity oscillator, and two controllers are proposed
for the variable order differential equations that track an arbitrary reference function.
The work \cite{Pedro} investigates the drag force acting on a particle due to the oscillatory
flow of a viscous fluid. The drag force is determined using the variable order fractional calculus,
where the order of derivative vary according to the dynamics of the flow.
In \cite{Ramirez:Coimbra2} a variable order differential equation
for a particle in a quiescent viscous liquid is developed.
For more on the application of variable order fractional operators
to the modeling of dynamic systems, we refer the reader to the recent
review article \cite{Ramirez:Coimbra1}.

In this note we develop the fractional calculus of variations
via a variable order approach. The~fractional variational calculus was born in 1996-1997
with the works of Riewe in mechanics \cite{CD:Riewe:1996,CD:Riewe:1997},
and is now under strong current research (see
\cite{MyID:182,DerInt,MyID:152,MyID:179,Cresson,gastao,fred:tor,MR1917624,MR1966935,Klimek02,%
Klimek05,comBasia,comDorota,MyID:181,MyID:207,MyID:226,MyID:227} and references therein).
However, to the best of the authors knowledge,
results for the variable order case are a rarity
and reduce to those in \cite{Atanackovic}.
Motivated by the advancements of \cite{DerInt,MyID:207},
and in contrast with \cite{Atanackovic},
we consider here fractional problems of the calculus
of variations where the Lagrangian depends
on classical integer order derivatives and
both on variable order fractional derivatives and integrals.

The paper is organized as follows. In \S\ref{sec:prelim} a brief review
to the variable order fractional calculus is given.
Our results are then formulated and proved in \S\ref{sec:MR}:
we show the boundedness of the variable order
Riemann--Liouville fractional integral in the space
$L_1[a,b]$ (Theorem~\ref{thm:1});
integration by parts formulas for variable order fractional operators
(Theorems~\ref{thm:IntByParts1} and \ref{thm:IntByParts2});
and a necessary optimality condition for our general fundamental problem
of the variable order fractional variational calculus (Theorem~\ref{thm:ELCaputo}).
Two illustrative examples are discussed in \S\ref{sec:ex}.


\section{Preliminaries}
\label{sec:prelim}

The reference book for fractional analysis and its applications is \cite{Samko:book}.
Here we recall the necessary definitions for the variable order
fractional calculus (see, \textrm{e.g.}, \cite{Lorentzo}).

\begin{defn}[Left and right Riemann--Liouville integrals of variable order]
Let $0<\alpha(t,\tau)<1$ for all $t, \tau \in [a,b]$ and $f\in L_1 [a,b]$. Then,
\begin{equation*}
{_{a}}\textsl{I}^{\alpha(\cdot,\cdot)}_{t}f(t)
= \int\limits_a^t\frac{1}{\Gamma(\alpha(t,\tau))}
(t-\tau)^{\alpha(t,\tau)-1}f(\tau)d\tau \quad (t>a)
\end{equation*}
is called the left Riemann--Liouville integral
of variable fractional order $\alpha(\cdot,\cdot)$, while
\begin{equation*}
{_{t}}\textsl{I}^{\alpha(\cdot,\cdot)}_{b}f(t)
=\int\limits_t^b \frac{1}{\Gamma(\alpha(\tau,t))}
(\tau-t)^{\alpha(\tau,t)-1}f(\tau)d\tau \quad (t<b)
\end{equation*}
denotes the right Riemann--Liouville integral
of variable fractional order $\alpha(\cdot,\cdot)$.
\end{defn}

\begin{defn}[Left and right Riemann--Liouville derivatives of variable order]
Let $0<\alpha(t,\tau)<1$ for all $t, \tau \in [a,b]$.
If ${_{a}}\textsl{I}^{1-\alpha(\cdot,\cdot)}_{t}f \in AC[a,b]$,
then the left Riemann--Liouville derivative of variable fractional
order $\alpha(\cdot,\cdot)$ is defined by
\begin{equation*}
{_{a}}\textsl{D}^{\alpha(\cdot,\cdot)}_{t} f(t)
= \frac{d}{dt} {_{a}}\textsl{I}^{1-\alpha(\cdot,\cdot)}_{t} f(t)
=\frac{d}{dt}\int\limits_a^t
\frac{1}{\Gamma(1-\alpha(t,\tau))}(t-\tau)^{-\alpha(t,\tau)} f(\tau)d\tau \quad (t>a)
\end{equation*}
while the right Riemann--Liouville derivative of variable fractional order
$\alpha(\cdot,\cdot)$ is defined for functions $f$ such that
${_{t}}\textsl{I}^{1-\alpha(\cdot,\cdot)}_{b}f\in AC[a,b]$ by
\begin{equation*}
{_{t}}\textsl{D}^{\alpha(\cdot,\cdot)}_{b}f(t)
= -\frac{d}{dt} {_{t}}\textsl{I}^{1-\alpha(\cdot,\cdot)}_{b}f(t)
=\frac{d}{dt}\int\limits_t^b
\frac{-1}{\Gamma(1-\alpha(\tau,t))}(\tau-t)^{-\alpha(\tau,t)}f(\tau)d\tau \quad (t<b).
\end{equation*}
\end{defn}

\begin{defn}[Left and right Caputo derivatives of variable fractional order]
\label{definition:Caputo}
Let $0<\alpha(t,\tau)<1$ for all $t, \tau \in [a,b]$.
If $f\in AC[a,b]$, then the left Caputo derivative
of variable fractional order $\alpha(\cdot,\cdot)$ is defined by
\begin{equation*}
{^{C}_{a}}\textsl{D}^{\alpha(\cdot,\cdot)}_{t}f(t)
=\int\limits_a^t
\frac{1}{\Gamma(1-\alpha(t,\tau))}(t-\tau)^{-\alpha(t,\tau)}\frac{d}{d\tau}f(\tau)d\tau \quad (t>a)
\end{equation*}
while the right Caputo derivative of variable fractional order $\alpha(\cdot,\cdot)$ is given by
\begin{equation*}
{^{C}_{t}}\textsl{D}^{\alpha(\cdot,\cdot)}_{b}f(t)
=\int\limits_t^b\frac{-1}{\Gamma(1-\alpha(\tau,t))}
(\tau-t)^{-\alpha(\tau,t)}\frac{d}{d\tau}f(\tau)d\tau \quad (t<b).
\end{equation*}
\end{defn}

The following result will be useful in the proof
of Theorems~\ref{thm:1} and \ref{thm:IntByParts1}.

\begin{thm}[\textrm{cf.} \cite{Ivady}]
Let $x \in [0,1]$. The Gamma function satisfies the inequalities
\begin{equation}
\label{eq:GammaIneq}
\frac{x^2+1}{x+1}\leq\Gamma(x+1)\leq\frac{x^2+2}{x+2}.
\end{equation}
\end{thm}


\section{Main results}
\label{sec:MR}

In \S\ref{sub:sec:bound} we prove that the variable fractional order
Riemann--Liouville integral ${_{a}}\textsl{I}^{\alpha(\cdot,\cdot)}_{t}$
is bounded on the space $L_1[a,b]$ (Theorem~\ref{thm:1}).
In \S\ref{sub:sec:ibp} we obtain a formula of integration by parts
for Riemann--Liouville integrals of variable order
(Theorem~\ref{thm:IntByParts1}) and a formula
of integration by parts for derivatives of variable fractional order
(Theorem~\ref{thm:IntByParts2}). Finally, in \S\ref{sub:sec:fvp}
we use the obtained formulas of integration by parts to derive
a necessary optimality condition for a general problem
of the calculus of variations  involving
variable order fractional operators (Theorem~\ref{thm:ELCaputo}).


\subsection{Boundedness}
\label{sub:sec:bound}

The following result allow us to consider
problems of the calculus of variations with a Lagrangian depending
on left Riemann--Liouville integrals of variable order.

\begin{thm}
\label{thm:1}
Let $\frac{1}{n}<\alpha(t,\tau)<1$ for all $t, \tau \in [a,b]$
and a certain $n\in\mathbb{N}$ greater or equal than two.
The Riemann--Liouville integral
${_{a}}\textsl{I}^{\alpha(\cdot,\cdot)}_{t}:L_1[a,b]\rightarrow L_1[a,b]$
of variable fractional order $\alpha(\cdot,\cdot)$
is a linear and bounded operator.
\end{thm}

\begin{proof}
The operator is obviously linear.
Let $\frac{1}{n}<\alpha(t,\tau)<1$, $f\in L_1 [a,b]$, and define
\[
F(\tau,t):=
\begin{cases}
\left|\frac{1}{\Gamma (\alpha(t,\tau))}(t-\tau)^{\alpha(t,\tau)-1}\right|
\cdot \left|f(\tau)\right| & \mbox{if $\tau < t$};\\
0 & \mbox{if $\tau \geq t$}
\end{cases}
\]
for all $(\tau,t)\in \Delta=[a,b]\times [a,b]$.
Since $\frac{1}{n}<\alpha(t,\tau)<1$, then
\begin{enumerate}

\item for $\tau+1\leq t$ we have $\ln(t-\tau)\geq 0$
and $(t-\tau)^{\alpha(t,\tau)-1}<1$;

\item for $\tau < t <\tau+1$ we have $\ln(t-\tau)<0$ and
$(t-\tau)^{\alpha(t,\tau)-1}<(t-\tau)^{\frac{1}{n}-1}$.
\end{enumerate}
Therefore,
\begin{equation*}
\begin{split}
\int_a^b &\left(\int_a^b F(\tau,t)dt\right)d\tau
=\int_a^b \left|f(\tau)\right|\left(\int_\tau^b\left|\frac{1}{\Gamma (\alpha(t,\tau))}
(t-\tau)^{\alpha(t,\tau)-1}\right|dt\right)d\tau\\
&=\int_a^b \left|f(\tau)\right|\left(\int_\tau^{\tau+1}\left|
\frac{(t-\tau)^{\alpha(t,\tau)-1}}{\Gamma(\alpha(t,\tau))}
\right|dt+\int_{\tau+1}^b\left|
\frac{(t-\tau)^{\alpha(t,\tau)-1}}{\Gamma(\alpha(t,\tau))}\right|
dt\right)d\tau\\
&<\int_a^b \left|f(\tau)\right|\left(\int_\tau^{\tau+1}\left|
\frac{(t-\tau)^{\frac{1}{n}-1}}{\Gamma(\alpha(t,\tau))}\right|dt
+\int_{\tau+1}^b\left|\frac{1}{\Gamma(\alpha(t,\tau))}\right|dt\right)d\tau.
\end{split}
\end{equation*}
Moreover, by inequality \eqref{eq:GammaIneq}, one has
\begin{equation*}
\begin{split}
&\int_a^b \left|f(\tau)\right|\left(\int_\tau^{\tau+1}\left|\frac{1}{\Gamma
(\alpha(t,\tau))}(t-\tau)^{\frac{1}{n}-1}\right|dt
+\int_{\tau+1}^b\left|\frac{1}{\Gamma (\alpha(t,\tau))}\right|dt\right)d\tau\\
&<\int_a^b \left|f(\tau)\right|\left(\int_\tau^{\tau+1}
\frac{\alpha^2(t,\tau)+\alpha(t,\tau)}{\alpha^2(t,\tau)+1}(t-\tau)^{\frac{1}{n}-1}dt
+\int_{\tau+1}^b\frac{\alpha^2(t,\tau)+\alpha(t,\tau)}{\alpha^2(t,\tau)+1}dt\right)d\tau\\
&<\int_a^b \left|f(\tau)\right|\left(\int_\tau^{\tau+1}
(t-\tau)^{\frac{1}{n}-1}dt+b-\tau-1\right)d\tau
=\int_a^b\left|f(\tau)\right|(n+b-\tau-1)d\tau\\
&<\left(n+ b - a \right)  \left\|f\right\|<\infty.
\end{split}
\end{equation*}
It follows from Fubini's theorem that $F$
is integrable in the square $\Delta$ and
\begin{equation*}
\begin{split}
\left\|{_{a}}\textsl{I}^{\alpha(\cdot,\cdot)}_{t} f\right\|
&=\int_a^b\left|
\int_a^t\frac{1}{\Gamma(\alpha(t,\tau))}(t-\tau)^{\alpha(t,\tau)-1}f(\tau)d\tau\right|dt\\
&\leq \int_a^b\left(\int_a^t\left|\frac{1}{\Gamma(\alpha(t,\tau))}
(t-\tau)^{\alpha(t,\tau)-1}f(\tau)\right|d\tau\right)dt\\
&=\int_a^b\left(\int_a^b F(\tau,t)d\tau\right)dt\\
&<\left(n+ b - a \right)\left\|f\right\|.
\end{split}
\end{equation*}
Therefore, ${_{a}}\textsl{I}^{\alpha(\cdot,\cdot)}_{t}:L_1[a,b]\rightarrow L_1[a,b]$
and $\left\|{_{a}}\textsl{I}^{\alpha(\cdot,\cdot)}_{t} \right\|< n+ b - a$.
\end{proof}


\subsection{Integration by parts formulas}
\label{sub:sec:ibp}

The integration by parts formulas we now obtain have an important role
in the proof of the generalized Euler--Lagrange equation \eqref{eq:eqELCaputo}.
We note that in Theorem~\ref{thm:IntByParts1} the left-hand side of \eqref{eq:ibp:1}
involves a left integral of variable order
while on the right-hand side it appears a right integral.

\begin{thm}
\label{thm:IntByParts1}
Let $\frac{1}{n}<\alpha(t,\tau)<1$ for all $t, \tau \in [a,b]$
and a certain $n\in\mathbb{N}$ greater or equal than two,
and $f,g\in C\left([a,b];\mathbb{R}\right)$. Then,
\begin{equation}
\label{eq:ibp:1}
\int_a^b g(t){_{a}}\textsl{I}^{\alpha(\cdot,\cdot)}_{t} f(t)dt
= \int_a^b f(t){_{t}}\textsl{I}^{\alpha(\cdot,\cdot)}_{b} g(t)dt.
\end{equation}
\end{thm}

\begin{proof}
Define
\[
F(\tau,t):= \left\{
\begin{array}{ll}
\left|\frac{1}{\Gamma(\alpha(t,\tau))}(t-\tau)^{\alpha(t,\tau)-1}
g(t)f(\tau)\right| & \mbox{if $\tau < t$};\\
0 & \mbox{if $\tau \geq t$}
\end{array} \right.
\]
for all $(\tau,t)\in \Delta = [a,b]\times [a,b]$.
Since $f$ and $g$ are continuous functions on $[a,b]$,
they are bounded on $[a,b]$, \textrm{i.e.},
there exist $C_1,C_2>0$ such that $\left|g(t)\right|\leq C_1$
and $\left|f(t)\right|\leq C_2$, $t\in[a,b]$. Therefore,
\begin{equation*}
\begin{split}
\int_a^b\left(\int_a^b F(\tau,t)d\tau\right)dt
&=\int_a^b \left(\int_a^t\left|\frac{1}{\Gamma(\alpha(t,\tau))}(t-\tau)^{\alpha(t,\tau)-1}
g(t)f(\tau)\right|d\tau\right)dt\\
&\leq C_1 C_2\int_a^b \left(\int_a^t\left|\frac{1}{\Gamma
\left(\alpha(t,\tau)\right)}(t-\tau)^{\alpha(t,\tau)-1}\right|d\tau\right)dt\\
&=C_1 C_2 \int_a^b \left(\int_a^t \frac{1}{\Gamma
(\alpha(t,\tau))} (t-\tau)^{\alpha(t,\tau)-1}d\tau\right)dt .
\end{split}
\end{equation*}
Since $\frac{1}{n}<\alpha(t,\tau)<1$, then
\begin{enumerate}

\item for $1\leq t-\tau $ we have $\ln(t-\tau)\geq 0$ and
$(t-\tau)^{\alpha(t,\tau)-1}<1$;

\item for $1>t-\tau$ we have $\ln(t-\tau)<0$ and
$(t-\tau)^{\alpha(t,\tau)-1}<(t-\tau)^{\frac{1}{n}-1}$.

\end{enumerate}
Therefore,
\begin{multline*}
C_1 C_2\int_a^b  \left(\int_a^t \frac{1}{\Gamma
\left(\alpha(t,\tau)\right)} (t-\tau)^{\alpha(t,\tau)-1}d\tau\right) dt\\
<C_1 C_2\int_a^b\left(\int_a^{t-1}\frac{1}{\Gamma
\left(\alpha(t,\tau)\right)} d\tau
+\int_{t-1}^t \frac{1}{\Gamma\left(\alpha(t,\tau)\right)}
(t-\tau)^{\frac{1}{n}-1} d\tau\right)dt.
\end{multline*}
Moreover, by \eqref{eq:GammaIneq}, one has
\begin{equation*}
\begin{split}
& C_1 C_2 \int_a^b\left(\int_a^{t-1}\frac{1}{\Gamma
\left(\alpha(t,\tau)\right)} d\tau
+\int_{t-1}^t \frac{1}{\Gamma
\left(\alpha(t,\tau)\right)}(t-\tau)^{\frac{1}{n}-1} d\tau\right)dt\\
&\leq C_1 C_2\int_a^b\left(\int_a^{t-1}\frac{\alpha^2(t,\tau)
+\alpha(t,\tau)}{\alpha^2(t,\tau)+1} d\tau
+\int_{t-1}^t \frac{\alpha^2(t,\tau)+\alpha(t,\tau)}{\alpha^2(t,\tau)
+1}(t-\tau)^{\frac{1}{n}-1} d\tau\right)dt\\
&< C_1 C_2\int_a^b\left(\int_a^{t-1}  d\tau
+\int_{t-1}^t (t-\tau)^{\frac{1}{n}-1} d\tau\right)dt\\
&=C_1 C_2(b-a)\left(\frac{b+a}{2}-1+n-a\right)<\infty.
\end{split}
\end{equation*}
Hence, one can use the Fubini theorem to change the order
of integration:
\begin{multline*}
\int_a^b g(t){_{a}}\textsl{I}^{\alpha(\cdot,\cdot)}_{t} f(t)dt
=\int_a^b \left(\int_a^t g(t)f(\tau)\frac{1}{\Gamma
\left(\alpha(t,\tau)\right)}(t-\tau)^{\alpha(t,\tau)-1}d\tau\right)dt\\
=\int_a^b \left(\int_\tau^b g(t)f(\tau)\frac{1}{\Gamma
\left(\alpha(t,\tau)\right)}(t-\tau)^{\alpha(t,\tau)-1}dt\right)d\tau
=\int_a^b f(\tau){_{\tau}}\textsl{I}^{\alpha(\cdot,\cdot)}_{b} g(\tau)d\tau.
\end{multline*}
\end{proof}

In our second formula \eqref{eq:intPart:Der} of fractional integration by parts,
the left-hand side contains a left Caputo derivative
of variable fractional order $\alpha(\cdot,\cdot)$, while on the
right-hand side it appears a right Riemann--Liouville integral
of variable order $1-\alpha(\cdot,\cdot)$ and a
right Riemann--Liouville derivative of variable order $\alpha(\cdot,\cdot)$.

\begin{thm}
\label{thm:IntByParts2}
Let $0 <\alpha(t,\tau)<1-\frac{1}{n}$ for all $t,\tau \in [a,b]$
and a certain $n\in\mathbb{N}$ greater or equal than two.
If $f\in C^1\left([a,b];\mathbb{R}\right)$,
$g\in C\left([a,b];\mathbb{R}\right)$,
and ${_{t}}\textsl{I}^{1-\alpha(\cdot,\cdot)}_{b} g \in AC[a,b]$, then
\begin{equation}
\label{eq:intPart:Der}
\int_a^b g(t){^{C}_{a}}\textsl{D}^{\alpha(\cdot,\cdot)}_{t}f(t)dt
=\left.f(t){_{t}}\textsl{I}^{1-\alpha(\cdot,\cdot)}_{b} g(t)\right|_a^b
+\int_a^b f(t){_{t}}\textsl{D}^{\alpha(\cdot,\cdot)}_{b}g(t)dt.
\end{equation}
\end{thm}

\begin{proof}
By Definition~\ref{definition:Caputo}, it follows that
${^{C}_{a}}\textsl{D}^{\alpha(\cdot,\cdot)}_{t}f(t)
={_{a}}\textsl{I}^{1-\alpha(\cdot,\cdot)}_{t}\frac{d}{dt} f(t)$.
Applying Theorem~\ref{thm:IntByParts1} and integration
by parts for classical (integer) derivatives, we obtain
\begin{equation*}
\begin{split}
\int_a^b g(t){^{C}_{a}}\textsl{D}^{\alpha(\cdot,\cdot)}_{t}f(t)dt
&=\int_a^b g(t){_{a}}\textsl{I}^{1-\alpha(\cdot,\cdot)}_{t}\frac{d}{dt}f(t)dt
=\int_a^b \frac{d}{dt}f(t){_{t}}\textsl{I}^{1-\alpha(\cdot,\cdot)}_{b}g(t)dt\\
&=\left.f(t){_{t}}\textsl{I}^{1-\alpha(\cdot,\cdot)}_{b} g(t)\right|_a^b
-\int_a^b f(t)\frac{d}{dt}{_{t}}\textsl{I}^{1-\alpha(\cdot,\cdot)}_{b}g(t)dt\\
&=\left.f(t){_{t}}\textsl{I}^{1-\alpha(\cdot,\cdot)}_{b} g(t)\right|_a^b
+\int_a^b f(t){_{t}}\textsl{D}^{\alpha(\cdot,\cdot)}_{b}g(t)dt.
\end{split}
\end{equation*}
\end{proof}


\subsection{A fundamental variational problem of variable fractional order}
\label{sub:sec:fvp}

We consider the problem of extremizing (minimizing or maximizing) a functional
\begin{equation}
\label{eq:31}
\mathcal{J}[y(\cdot)]=\int\limits_a^b
F\left(t,y(t),y'(t),{^{C}_{a}}\textsl{D}^{\alpha(\cdot,\cdot)}_{t}y(t),
{_{a}}\textsl{I}^{\beta(\cdot,\cdot)}_{t}y(t)\right)dt
\end{equation}
subject to boundary conditions
\begin{equation}
\label{eq:32}
y(a)=y_a, \quad y(b)=y_b,
\end{equation}
where $\alpha, \beta : [a,b]\times [a,b] \rightarrow \mathbb{R}$
are given functions taking values in $(0,1-\frac{1}{n})$
and $(\frac{1}{n}, 1)$, respectively,
with $n\in\mathbb{N}$ greater or equal than two.
For simplicity of notation, we introduce the operator
$\left\{ \cdot, \cdot, \cdot \right\}$ defined by
\begin{equation*}
\left\{y,\alpha,\beta\right\}(t)
=\left(t,y(t),y'(t),{^{C}_{a}}\textsl{D}^{\alpha(\cdot,\cdot)}_{t}y(t),
{_{a}}\textsl{I}^{\beta(\cdot,\cdot)}_{t}y(t)\right).
\end{equation*}
We assume that $F\in C^1\left([a,b]\times\mathbb{R}^4;\mathbb{R}\right)$;
$t \mapsto \partial_4 F\left\{y,\alpha,\beta\right\}(t)$ is continuous and
has absolutely continuous integral ${_{t}}\textsl{I}^{1-\alpha(\cdot,\cdot)}_{b}$
and continuous derivative ${_{t}}\textsl{D}^{\alpha(\cdot,\cdot)}_{b}$;
$t \mapsto \partial_5 F\left\{y,\alpha,\beta\right\}(t)$ is continuous and
has continuous variable order fractional integral ${_{t}}\textsl{I}^{\beta(\cdot,\cdot)}_{b}$;
and $t \mapsto \partial_3 F\left\{y,\alpha,\beta\right\}(t)$
has continuous usual derivative $\frac{d}{dt}$.

\begin{defn}
A continuously differentiable function $y\in C^1\left([a,b];\mathbb{R}\right)$ is said to be
admissible for the variational problem \eqref{eq:31}--\eqref{eq:32},
if ${^{C}_{a}}\textsl{D}^{\alpha(\cdot,\cdot)}_{t} y$
and ${_{a}}\textsl{I}^{\beta(\cdot,\cdot)}_{t} y$ exist and are continuous on the
interval $[a,b]$, and $y$ satisfies the given boundary conditions \eqref{eq:32}.
\end{defn}

\begin{thm}
\label{thm:ELCaputo}
Let $y$ be a solution to problem \eqref{eq:31}--\eqref{eq:32}.
Then, $y$ satisfies the generalized Euler--Lagrange equation
\begin{multline}
\label{eq:eqELCaputo}
\partial_2 F\left\{y,\alpha,\beta\right\}(t)
-\frac{d}{dt}\partial_3 F\left\{y,\alpha,\beta\right\}(t)
+{_{t}}\textsl{I}^{\beta(\cdot,\cdot)}_{b}\partial_5 F\left\{y,\alpha,\beta\right\}(t)\\
+{_{t}}\textsl{D}^{\alpha(\cdot,\cdot)}_{b}\partial_4 F\left\{y,\alpha,\beta\right\}(t)=0
\end{multline}
for all $t\in[a,b]$.
\end{thm}

\begin{proof}
Suppose that $y$ is an extremizer of $\mathcal{J}$.
Consider the value of $\mathcal{J}$ at a nearby function
$\hat{y}(t)=y(t)+\varepsilon\eta(t)$,
where $\varepsilon\in\mathbb{R}$ is a small parameter
and $\eta\in C^1\left([a,b];\mathbb{R}\right)$ is an arbitrary
function satisfying $\eta(a)=\eta(b)=0$ and such that
${^{C}_{a}}\textsl{D}^{\alpha(\cdot,\cdot)}_{t} \hat{y}(t)$
and ${_{a}}\textsl{I}^{\beta(\cdot,\cdot)}_{t} \hat{y}(t)$ are continuous. Let
$$
J(\varepsilon)=\mathcal{J}[\hat{y}(\cdot)]
=\int_a^b F\left\{\hat{y},\alpha,\beta\right\}(t) dt.
$$
A necessary condition for $y$ to be an extremizer is given by
\begin{multline}
\label{eq:33}
\left.\frac{d J}{d\varepsilon}\right|_{\varepsilon=0}=0
\Leftrightarrow\int\limits_a^b
\Biggl(\partial_2 F\left\{y,\alpha,\beta\right\}(t) \cdot \eta(t)
+\partial_3 F \left\{y,\alpha,\beta\right\}(t)\frac{d}{dt}\eta(t)\\
+ \partial_4 F\left\{y,\alpha,\beta\right\}(t)
{^{C}_{a}}\textsl{D}^{\alpha(\cdot,\cdot)}_{t}\eta(t)
+\partial_5 F \left\{y,\alpha,\beta\right\}(t)
\cdot {_{a}}\textsl{I}^{\beta(\cdot,\cdot)}_{t} \eta(t)\Biggr)dt = 0.
\end{multline}
Using the classical and the generalized fractional integration by parts formulas
of Theorems~\ref{thm:IntByParts1} and \ref{thm:IntByParts2}, we obtain
\begin{equation*}
\int_a^b\partial_3F\frac{d\eta}{dt}dt=\left.\partial_3F\eta\right|_a^b
-\int_a^b\left(\eta\frac{d}{dt}\partial_3F\right)dt,
\end{equation*}
\begin{equation*}
\int\limits_a^b\partial_4 F {^{C}_{a}}\textsl{D}^{\alpha(\cdot,\cdot)}_{t}\eta dt
=\left.\eta {_{t}}\textsl{I}^{1-\alpha(\cdot,\cdot)}_{b}\partial_4 F\right|_a^b
+\int\limits_a^b\eta{_{t}}\textsl{D}^{\alpha(\cdot,\cdot)}_{b}\partial_4 F dt,
\end{equation*}
and
\begin{equation*}
\int\limits_a^b\partial_5 F {_{a}}\textsl{I}^{\beta(\cdot,\cdot)}_{t}\eta dt
=\int\limits_a^b
\eta {_{t}}\textsl{I}^{\beta(\cdot,\cdot)}_{b}\partial_5 F dt.
\end{equation*}
Because $\eta(a)=\eta(b)=0$, \eqref{eq:33} simplifies to
\begin{multline*}
\int_a^b\eta(t)\Biggl(\partial_2 F \left\{y,\alpha,\beta\right\}(t)
-\frac{d}{dt}\partial_3 F\left\{y,\alpha,\beta\right\}(t)
+{_{t}}\textsl{D}^{\alpha(\cdot,\cdot)}_{b}\partial_4 F\left\{y,\alpha,\beta\right\}(t)\\
+{_{t}}\textsl{I}^{\beta(\cdot,\cdot)}_{b}\partial_5 F\left\{y,\alpha,\beta\right\}(t)\Biggr)dt=0.
\end{multline*}
One obtains \eqref{eq:eqELCaputo}
by the fundamental lemma of the calculus of variations
(see, \textrm{e.g.}, \cite{book:vanBrunt}).
\end{proof}


\section{Illustrative examples}
\label{sec:ex}

Let $\beta(t,\tau) = \beta(t)$ be a function
depending only on variable $t$,
$\frac{1}{n}<\beta(t)<1$ for all $t \in [a,b]$
and a certain $n\in\mathbb{N}$ greater or equal than two,
and $\gamma>-1$. We make use of the identity
\begin{equation}
\label{eq:VOFDPower}
{_{a}}\textsl{I}^{\beta(\cdot)}_{t} (t-a)^{\gamma}
=\frac{\Gamma(\gamma+1)(t-a)^{\gamma+\beta(t)}}{\Gamma(\gamma+\beta(t)+1)}
\end{equation}
that one can find in the Samko and Ross paper \cite{Samko:Ross2}.

\begin{ex}
\label{ex:VOFVCExample}
Let $\mathcal{J}$ be the functional defined by
\begin{equation*}
\mathcal{J}[y(\cdot)]=\int\limits_a^b \sqrt{1+\frac{\Gamma(\beta(t)+3)}{2\Gamma(3)(t-a)^{2
+\beta(t)}}\left({_{a}}\textsl{I}^{\beta(\cdot)}_{t}y(t))\right)^2
-{_{a}}\textsl{I}^{\beta(\cdot)}_{t} y(t)} \, dt
\end{equation*}
with endpoint conditions $y(a)=0$ and $y(b)=(b-a)^2$. If $y$ is an extremizer
for $\mathcal{J}$, then the necessary optimality condition
of Theorem~\ref{thm:ELCaputo} gives
\begin{equation}
\label{eq:ELExample}
{_{t}}\textsl{I}^{\beta(\cdot)}_{b}\left(\frac{\frac{\Gamma(\beta(t)
+3)}{\Gamma(3)(t-a)^{2+\beta(t)}}{_{a}}\textsl{I}^{\beta(\cdot)}_{t} y(t)
-1}{2\sqrt{1+\frac{\Gamma(\beta(t)+3)}{2\Gamma(3)(t-a)^{2
+\beta(t)}}\left({_{a}}\textsl{I}^{\beta(\cdot)}_{t}y(t)\right)^2
-{_{a}}\textsl{I}^{\beta(\cdot)}_{t}y(t)}}\right)=0.
\end{equation}
By identity \eqref{eq:VOFDPower}, function
\begin{equation}
\label{eq:solution}
y(t)=(t-a)^{2}
\end{equation}
is a solution to the variable order fractional differential equation
\begin{equation*}
{_{a}}\textsl{I}^{\beta(\cdot)}_{t}y(t)
=\frac{\Gamma(3)(t-a)^{2+\beta(t)}}{\Gamma(\beta(t)+3)}.
\end{equation*}
Therefore, \eqref{eq:solution} is a solution
to the Euler--Lagrange equation \eqref{eq:ELExample}
and an extremal.
\end{ex}

In the next example $\alpha(t,\tau)$ is a function taking values in the set
$(0,1-\frac{1}{n})$.

\begin{ex}
Consider the following problem:
\begin{gather*}
\mathcal{J}[y(\cdot)]=\int\limits_a^b
\left({^{C}_{a}}\textsl{D}^{\alpha(\cdot,\cdot)}_{t}y(t)\right)^2
+ \left({_{a}}\textsl{I}^{\beta(\cdot)}_{t}y(t)
- \frac{\xi(t-\tau)^{\beta(t)}}{\Gamma(\beta(t)+1)}\right)^2 dt \longrightarrow \min,\\
y(a) = \xi, \quad y(b) = \xi,
\end{gather*}
for a given real $\xi$. Because $\mathcal{J}[y(\cdot)] \geq 0$ for any function $y$
and $\mathcal{J}[\tilde{y}(\cdot)] = 0$ for the admissible function $\tilde{y} = \xi$
(use relation \eqref{eq:VOFDPower} for $\gamma=0$, the linearity of operator
${_{a}}\textsl{I}^{\beta(\cdot)}_{t}$, and the definition of left Caputo derivative
of a variable fractional order), we conclude that $\tilde{y}$
is the global minimizer to the problem. It is straightforward
to check that $\tilde{y}$ satisfies our variable order fractional
Euler--Lagrange equation \eqref{eq:eqELCaputo}.
\end{ex}


\subsection*{Acknowledgments}

Work supported by {\it FEDER} funds through
{\it COMPETE} --- Operational Programme Factors of Competitiveness
(``Programa Operacional Factores de Competitividade'')
and by Portuguese funds through the
{\it Center for Research and Development
in Mathematics and Applications} (University of Aveiro)
and the Portuguese Foundation for Science and Technology
(``FCT --- Funda\c{c}\~{a}o para a Ci\^{e}ncia e a Tecnologia''),
within project PEst-C/MAT/UI4106/2011
with COMPETE number FCOMP-01-0124-FEDER-022690.
Odzijewicz was also supported by FCT through the Ph.D. fellowship
SFRH/BD/33865/2009; Malinowska by Bia{\l}ystok
University of Technology grant S/WI/02/2011;
and Torres by FCT through the project PTDC/MAT/113470/2009.




\begin{thebibliography}{xx}

\bibitem{alexandre}
A. Almeida, S. Samko,
\textit{Fractional and hypersingular operators
in variable exponent spaces on metric measure spaces.}
Mediterr. J. Math. \textbf{6} (2009), 215--232.

\bibitem{MyID:182}
R. Almeida, A. B. Malinowska, D. F. M. Torres,
\textit{A fractional calculus of variations for multiple
integrals with application to vibrating string.}
J. Math. Phys. \textbf{51} (2010), 033503, 12~pp.
{\tt arXiv:1001.2722}

\bibitem{DerInt}
R. Almeida, D. F. M. Torres,
\textit{Calculus of variations with fractional
derivatives and fractional integrals.}
Appl. Math. Lett. \textbf{22} (2009), 1816--1820.
{\tt arXiv:0907.1024}

\bibitem{Atanackovic}
T. M. Atanackovic, S. Pilipovic,
\textit{Hamilton's principle with variable order fractional derivatives.}
Fract. Calc. Appl. Anal. \textbf{14} (2011), 94--109.

\bibitem{MyID:152}
N. R. O. Bastos, R. A. C. Ferreira, D. F. M. Torres,
\textit{Necessary optimality conditions for fractional
difference problems of the calculus of variations.}
Discrete Contin. Dyn. Syst. \textbf{29} (2011), 417--437.
{\tt arXiv:1007.0594}

\bibitem{MyID:179}
N. R. O. Bastos, R. A. C. Ferreira, D. F. M. Torres,
\textit{Discrete-time fractional variational problems.}
Signal Process. \textbf{91} (2011), 513--524.
{\tt arXiv:1005.0252}

\bibitem{Coimbra}
C. F. M. Coimbra,
\textit{Mechanics with variable-order differential operators.}
Ann. Phys. \textbf{12} (2003), 692--703.

\bibitem{Cresson}
J. Cresson,
\textit{Fractional embedding of differential operators and Lagrangian systems.}
J. Math. Phys. \textbf{48} (2007), 033504, 34~pp.
{\tt arXiv:math/0605752}

\bibitem{Diaz:Coimbra}
G. Diaz, C. F. M. Coimbra,
\textit{Nonlinear dynamics and control of a variable order
oscillator with application to the van der Pol equation.}
Nonlinear Dynam. \textbf{56} (2009), 145--157.

\bibitem{gastao}
G. S. F. Frederico, D. F. M. Torres,
\textit{Fractional conservation laws in optimal control theory.}
Nonlinear Dynam. \textbf{53} (2008), 215--222.
{\tt arXiv:0711.0609}

\bibitem{fred:tor}
G. S. F. Frederico, D. F. M. Torres,
\textit{Fractional Noether's theorem in the Riesz-Caputo sense.}
Appl. Math. Comput. \textbf{217} (2010), 1023--1033.
{\tt arXiv:1001.4507}

\bibitem{Ivady}
P. Ivady,
\textit{A note on a gamma function inequality.}
J. Math. Inequal. \textbf{3} (2009), 227--236.

\bibitem{book:Kilbas}
A. A. Kilbas, H. M. Srivastava, J. J. Trujillo,
\textit{Theory and applications of fractional differential equations},
North-Holland Mathematics Studies, 204, Elsevier, Amsterdam, 2006.

\bibitem{MR1917624}
M. Klimek,
\textit{Fractional sequential mechanics---models with symmetric fractional derivative},
Czechoslovak J. Phys. {\bf 51} (2001), no.~12, 1348--1354.

\bibitem{MR1966935}
M. Klimek,
\textit{Lagrangean and Hamiltonian fractional sequential mechanics},
Czechoslovak J. Phys. {\bf 52} (2002), no.~11, 1247--1253.

\bibitem{Klimek02}
M. Klimek,
\textit{Stationarity-conservation laws for fractional
differential equations with variable coefficients.}
J. Phys. A \textbf{35} (2002), 6675--6693.

\bibitem{Klimek05}
M. Klimek,
\textit{Lagrangian fractional mechanics---a noncommutative approach.}
Czechoslovak J. Phys. \textbf{55} (2005), 1447--1453.

\bibitem{book:Klimek}
M. Klimek,
\textit{On solutions of linear fractional differential equations of a variational type.}
The Publishing Office of Czenstochowa University of Technology, Czestochowa, 2009.

\bibitem{Lorentzo}
C. F. Lorenzo, T. T. Hartley,
\textit{Variable order and distributed order fractional operators.}
Nonlinear Dynam. \textbf{29} (2002), 57--98.

\bibitem{comBasia}
A. B. Malinowska, D. F. M. Torres,
\textit{Generalized natural boundary conditions
for fractional variational problems
in terms of the Caputo derivative.}
Comput. Math. Appl. \textbf{59} (2010), 3110--3116.
{\tt arXiv:1002.3790}

\bibitem{comDorota}
D. Mozyrska, D. F. M. Torres,
\textit{Minimal modified energy control for fractional
linear control systems with the Caputo derivative.}
Carpathian J. Math. \textbf{26} (2010), 210--221.
{\tt arXiv:1004.3113}

\bibitem{MyID:181}
D. Mozyrska, D. F. M. Torres,
\textit{Modified optimal energy and initial memory
of fractional continuous-time linear systems.}
Signal Process. \textbf{91} (2011), 379--385.
{\tt arXiv:1007.3946}

\bibitem{MyID:207}
T. Odzijewicz, A. B. Malinowska, D. F. M. Torres,
\textit{Fractional variational calculus
with classical and combined Caputo derivatives.}
Nonlinear Anal. \textbf{75} (2012), 1507--1515.
{\tt arXiv:1101.2932}

\bibitem{MyID:226}
T. Odzijewicz, A. B. Malinowska, D. F. M. Torres,
\textit{Generalized fractional calculus with applications
to the calculus of variations},
Comput. Math. Appl. (2012), in press.
DOI: 10.1016/j.camwa.2012.01.073
{\tt arXiv:1201.5747}

\bibitem{MyID:227}
T. Odzijewicz, A. B. Malinowska, D. F. M. Torres,
\textit{Fractional calculus of variations in terms
of a generalized fractional integral with applications to physics},
Abstr. Appl. Anal. (2012), in press.
{\tt arXiv:1203.1961}

\bibitem{Pedro}
H. T. C. Pedro, M. H. Kobayashi, J. M. C. Pereira, C. F. M. Coimbra,
\textit{Variable order modeling of diffusive-convective effects
on the oscillatory flow past a sphere.}
J. Vib. Control \textbf{14} (2008), 1569--1672.

\bibitem{book:Podlubny}
I. Podlubny,
\textit{Fractional differential equations.}
Mathematics in Science and Engineering, 198,
Academic Press, San Diego, CA, 1999.

\bibitem{Ramirez:Coimbra1}
L. E. S. Ramirez, C. F. M. Coimbra,
\textit{On the selection and meaning
of variable order operators for dynamic modeling.}
Int. J. Differ. Equ. \textbf{2010} (2010), Art. ID 846107, 16~pp.

\bibitem{Ramirez:Coimbra2}
L. E. S. Ramirez and C. F. M. Coimbra,
\textit{On the variable order dynamics of the
nonlinear wake caused by a sedimenting particle.}
Phys. D \textbf{240} (2011), 1111--1118.

\bibitem{CD:Riewe:1996}
F. Riewe,
\textit{Nonconservative Lagrangian and Hamiltonian mechanics.}
Phys. Rev. E (3) \textbf{53} (1996), 1890--1899.

\bibitem{CD:Riewe:1997}
F. Riewe,
\textit{Mechanics with fractional derivatives.}
Phys. Rev. E (3) \textbf{55} (1997), 3581--3592.

\bibitem{Samko:Ross1}
B. Ross, S. G. Samko,
\textit{Fractional integration operator of a variable
order in the Holder spaces $H^{\lambda(x)}$.}
Internat. J. Math. Math. Sci. \textbf{18} (1995), 777--788.

\bibitem{Samko}
S. G. Samko,
\textit{Fractional integration and differentiation of variable order.}
Anal. Math. \textbf{21} (1995), 213--236.

\bibitem{Samko:book}
S. G. Samko, A. A. Kilbas, O. I. Marichev,
\textit{Fractional integrals and derivatives},
translated from the 1987 Russian original, Gordon and Breach, Yverdon, 1993.

\bibitem{Samko:Ross2}
S. G. Samko, B. Ross,
\textit{Integration and differentiation to a variable fractional order.}
Integral Transform. Spec. Funct. \textbf{1} (1993), 277--300.

\bibitem{TM:rec}
J. A. Tenreiro Machado, V. Kiryakova, F. Mainardi,
\textit{Recent history of fractional calculus.}
Commun. Nonlinear Sci. Numer. Simul. \textbf{16} (2011), 1140--1153.

\bibitem{book:vanBrunt}
B. van Brunt,
\textit{The calculus of variations.}
Universitext, Springer, New York, 2004.

\end{thebibliography}
\end{document}